\pdfoutput=1
\documentclass[]{article}
\usepackage{amssymb,amsmath,hyperref,verbatim,graphicx,amsthm}
\makeatletter \let\cl@chapter\relax \makeatother
\usepackage{cleveref, xurl}
\usepackage{multirow, longtable}
\usepackage{booktabs}
\usepackage{diagbox}
\usepackage{authblk}
\newtheorem{theorem}{Theorem}[section]
\newtheorem{proposition}[theorem]{Proposition}
\newtheorem{corollary}[theorem]{Corollary}
\newtheorem{lemma}[theorem]{Lemma}
\usepackage{pdflscape}
\usepackage{mathtools}

\usepackage{enumerate}
\usepackage[dvipsnames]{xcolor}
\usepackage{tikz}
\usepackage{pgfplots}
\pgfplotsset{width=10cm,compat=1.9}
\usepackage{amsmath, amssymb}
\usepackage{booktabs}
\usepackage[ruled,vlined,linesnumbered]{algorithm2e}
\usepackage{hyperref}
\usepackage{subcaption}
\hypersetup{
    colorlinks=true,
    linkcolor=blue,
    filecolor=magenta,      
    urlcolor=cyan,
}

\begin{document}

\newcommand{\red}{\color{red}}
\newcommand{\blue}{\color{blue}}
\newcommand{\green}{\color{green}}
\allowdisplaybreaks

\title{Parallelized Conflict Graph Cut Generation}
\author[]{Yongzheng Dai}
\author[]{Chen Chen}
\affil[]{ISE, The Ohio State University, Columbus, OH, USA}

\maketitle

\begin{abstract}
A conflict graph represents logical relations between binary variables, and effective use of the graph can significantly accelerate branch-and-cut solvers for mixed-integer programming (MIP). In this paper we develop efficient parallel conflict graph management: conflict detection; maximal clique generation; clique extension; and clique merging. We leverage parallel computing in order to intensify computational effort on the conflict graph, thereby generating a much larger pool of cutting planes than what can be practically achieved in serial. Computational experiments demonstrate that the expanded pool of cuts enabled by parallel computing lead to substantial reductions in total MIP solve time, especially for more challenging cases.
\end{abstract}

\section{Introduction}
Various progress and benchmarking reports from the literature \cite{bixby2007progress,berthold2012solving,achterberg2013mixed,koch2022progress} indicate the increasing importance of parallelization as well as preprocessing subroutines in solvers for mixed-integer programming (MIP). We build upon this progress by parallelizing conflict graph (CG) \cite{scip} management, a key subroutine in branch-and-cut solvers \cite{hoffman1993solving,atamturk2000conflict,achterberg2007conflict,brito2021preprocessing} that begins at the preprocessing stage and can be subsequently deployed throughout the search tree (see, e.g. \cite{witzig2017experiments}). Specifically, we adopt the CG definition of Achterberg \cite{scip} (referred to as a clique graph in SCIP \cite{berthold2012solving}), in contrast with e.g. \cite{achterberg2007conflict,witzig2017experiments}. This work focuses on modest levels of parallelism---say, $<100$ cores with shared memory---typical of personal computing setups now and in the near future; in contrast, work on massively parallel MIP (e.g. \cite{phillips2006massively,eckstein2015pebbl,shinano2016solving,koch2012could,shinano2018fiberscip,perumalla2021design}) address certain issues associated with distributed computing such as higher communication costs. Moreover, our work deals with general-purpose MIP; for instance, certain stochastic optimization problems have special problem structures amenable to specially tailored decomposition-based distributed schemes (see e.g. \cite{papavasiliou2014applying,ryan2016scenario,munguia2019parallel}). 

Perhaps the most closely related work on parallel MIP is by Gleixner, Gottwald, and Hoen \cite{gleixner2023papilo}, in which a host of preprocessing techniques are parallelized. Using 32 threads, their solver PaPILO is reported to reduce presolve time by almost 50\% in shifted geometric mean, with 5x speedup on preprocessing-intensive instances. Our method attains over 80\% reduction using 64 threads (see Figure~\ref{fig:parallel_comp}); we emphasize, however, that a direct comparison is misguided as both our works are complementary due to parallelization of different procedures. Likewise, our work complements other efforts such as parallel LP solvers (e.g. \cite{klabjan2000parallel,huangfu2018parallelizing,andersen1996implementation,hall2010towards}), and concurrent solves (see e.g. \cite{koch2022progress}).  Our approach is further differentiated from the literature as we do not solely accelerate existing serial algorithms---indeed, the  small fraction of total runtime from typical CG management suggests rather modest potential from this due to Amdahl's law \cite{gustafson1988reevaluating}--- but instead modify the serial CG procedures of Brito and Santos \cite{brito2021preprocessing} to generate more cuts. We observe empirically that our more intensive cut management scheme is only modestly above break-even in serial implementation (time spent in cut generation offsets roughly equal amounts of time in the branch-and-cut solver), but attains substantial overall time speedups when executed in parallel. 

\section{Parallel Conflict Graph Management}
The serial algorithm development in this section follows predominantly the recent work of Brito and Santos \cite{brito2021preprocessing} on CG management, which has been implemented in the COIN-OR branch-and-cut solver.

Furthermore, throughout the paper, in parallel computing analysis we assume $k$ threads with shared memory among cores.  Moreover, theoretical results herein hold irrespective of the particular PRAM configuration (EREW, CRCW, etc.) due to the lack of conflicting transactions in our algorithms.

Consider a mixed integer program (MIP) in the following generic form:
\begin{equation}\label{eq:MIP}
    \begin{aligned}
        \min\ &c^Tx\\
        \mbox{s.t.}\quad &Ax \circ b\\
        &\ell \leq x\leq u\\
        &x_j\in \mathbb{Z} \mbox{ for all } j \in \mathcal{I}
    \end{aligned}
\end{equation}
with parameters $A\in \mathbb{R}^{m\times n}$, $c\in \mathbb{R}^n$, $b\in \mathbb{R}^m$, $\ell \in (\mathbb{R}\cup \{-\infty\})^n$, and $u \in (\mathbb{R}\cup \{\infty\})^n$; variables $x\in \mathbb{R}^n$ with $x_j\in \mathbb{Z}$ for $j\in \mathcal{I} \subseteq \mathcal{N} = \{1,...,n\}$; and constraints with relations $\circ_i\in\{=,\leq,\geq\}$ for each row $i\in \mathcal{M} = \{1,...,m\}$. Furthermore, let $\mathcal{B} := \{j\in \mathcal{I}\ |\ \ell_j=0 \wedge u_j = 1\}$ be the set of indices for all binary variables.

A conflict graph over $\mathcal{B}$ has one node for each binary variable $x_j$, representing an assigned value $x_j=1$, and another node for the complement variable $\bar x_j :=1-x_j$ assigned to 1. An edge between the nodes represents a conflict, i.e., an infeasible assignment; for instance, one can automatically place an edge between each variable and its complement. Conflicts need to be inferred or detected from the problem formulation or else possibly during the branch-and-bound procedure (see, e.g.,\cite{witzig2021computational}). In this paper, we consider two types of constraints from which conflicts can be extracted: \emph{set packing} constraints and \emph{conflicting knapsack} constraints. A set packing constraint is defined as
\begin{equation}\label{eq:setpack}
    \sum_{j\in \mathcal{S}}x_j \leq 1,
\end{equation}
for some $\mathcal{S}\subseteq \mathcal{B}$. Since each variable in $\mathcal{S}$ has a conflict with all others, $\mathcal{S}$ forms a clique in the conflict graph.

A \emph{knapsack constraint} is defined as
\begin{equation}\label{eq:knapsack}
    \sum_{j\in \mathcal{B}}a_jx_j \leq b
\end{equation}
with $a_j \geq 0, j\in \mathcal{B}$. Suppose WLOG that $a_{(1)}, a_{(2)}$ are the two largest elements from coefficients $(a_j)_{j\in \mathcal{B}}$. Then if $a_{(1)}+a_{(2)} > b$, we call this a \emph{conflicting knapsack constraint} as a CG clique can be generated from $x_1,x_2$ (and possibly other variables); in the absence of this condition, no conflicts can be inferred from the knapsack.

Set packing and knapsack constraints can, in turn, be inferred from the general MIP formulation. For a given mixed integer constraint $A_{i\cdot}x \circ_i b_i$, if $\circ_i$ is $\leq$, we extract a pure binary variables constraint (PBC) denoted as $\mbox{PBC}(i)$:
\begin{equation}\label{eq:pure_binary}
    \sum_{j\in \mathcal{B}: A_{ij}>0}A_{ij}x_j - \sum_{j\in \mathcal{B}: A_{ij}<0}A_{ij}\bar{x}_j \leq b_i - \inf\{\sum_{j\not\in \mathcal{B}}A_{ij}x_j\} - \sum_{j\in \mathcal{B}: A_{ij}<0}A_{ij}.
\end{equation}
If $\inf\{\sum_{j \notin B} A_{ij} x_j\} = -\infty$, $\mbox{PBC}(i)$ does not infer any conflicts and can be ignored. Note that if $\circ_i$ is $\geq$, we can rewrite the constraint as $-A_{i\cdot}x \leq -b_i$, and if $\circ_i$ is $=$, we can split the constraint into two constraints $A_{i\cdot}x \leq b_i$ and $-A_{i\cdot}x \leq -b_i$. So we only consider constraints with the form of $A_{i\cdot}x \leq b_i$ in the remainder of this section.

The remainder of this section describes both serial and parallel algorithms for certain aspects of CG management, presented in order of execution during preprocessing.
Both serial and parallel algorithms are analyzed in terms of worst-case and average-case complexity.  

\subsection{Detecting Set Packing and Conflicting Knapsack Constraints}
\label{sec:pbc}
At the preprocessing stage, given a $\mbox{MIP} = (\mathcal{M}, \mathcal{N}, \mathcal{I}, A, b, c, \circ, \ell, u)$, we perform one-round detection of set packing constraints and conflicting knapsack constraints (subsequently used for conflict graph cut generation), which is described in Alg~\ref{Alg:simple_presolve}. This applies well-studied techniques from the CG literature (see e.g. Achterberg et al. \cite{achterberg2020presolve}). We replicate this standard component of presolve for implementation purposes in order to develop a standalone CG procedure for testing, as we do not have internal access to any MIP solvers.


In line 4 of Alg~\ref{Alg:simple_presolve}, if $A_{i\cdot}x \leq b_i$ is a set packing constraint, we call it an \emph{original} set packing (OSP) constraint and collect all such constraints in the set $\mathcal{S}_{osp}$ (line 6). These set packing constraints are removed from $\mathcal{M}$ as they are later replaced with strengthened versions (see Section~\ref{sec:cutmgmt}).

In line 10, if PBC(i) is a singleton $A_{ij} x_j \leq b_i$ for the binary variable $x_j$, then PBC(i) is either redundant or we can fix $x_j$ to $0$ or $1$, depending on the values of $A_{ij}$ and $b_i$ (e.g. domain propagation \cite{achterberg2020presolve}).

In line 11 of Alg~\ref{Alg:simple_presolve}, if $\mbox{PBC}(i)$ is a set packing constraint, we call it an \emph{inferred} set packing constraint (ISP) and collect all such constraints in the set $\mathcal{S}_{isp}$.

In line 14 we collect all (i.e. both original and PBC-inferred) conflicting knapsack constraints (CK) in the constraint set $\mathcal{S}_{ck}$. 

\begin{algorithm}[ht]
\SetAlgoLined
\LinesNumbered
\SetKwRepeat{Do}{do}{while}
\SetKwInput{Input}{Input}
\SetKwInput{Output}{Output}
\Input{$\mbox{MIP} = (\mathcal{M}, \mathcal{N}, \mathcal{I}, A, b, c, \circ, \ell, u)$}
\Output{MIP after applying Simple Presolve, $\mathcal{S}_{osp}$, $\mathcal{S}_{isp}$, $\mathcal{S}_{ck}$}
 Set $\mathcal{I}_{osp}$, $\mathcal{S}_{isp}$, $\mathcal{S}_{ck}:=\emptyset$\;
 Remove empty constraints and singletons from $\mbox{MIP}$ and conduct a one-round single-row bound strengthening to $\mbox{MIP}$\;
 \For{$i \in \mathcal{M}$}{
    \eIf{$A_{i\cdot}x \leq b_i$ is a set packing constraint}{
    Remove $i$ from $\mathcal{M}$\;
    Set $\mathcal{S}_{osp}:= \mathcal{S}_{osp}\cup\{A_{i\cdot}x\leq b_i\}$\;}
    {
        Rewrite $A_{i\cdot}x \leq b_i$ as $\mbox{PBC}(i)$\;
        \uIf{$\mathrm{PBC}(i)$ is a singleton}{Update $u_j$ and $\ell_j$ for $j$ such that $A_{ij} \neq 0$ and $j\in \mathcal{B}$\;}
        \uElseIf{$\mathrm{PBC}(i)$ is a set packing constraint}{Set $\mathcal{S}_{isp} := \mathcal{S}_{isp}\cup\{\mbox{PBC}(i)\}$\;}
        \Else
        {Set $\mathcal{S}_{ck}:=\mathcal{S}_{ck}\cup\{\mbox{PBC}(i)\}$\;}
    }
 }
 \textbf{return} MIP after applying Simple Presolve, $\mathcal{S}_{osp}$, $\mathcal{S}_{isp}$, $\mathcal{S}_{ck}$. 
 \caption{Detecting Set Packing and Conflicting Knapsack Constraints}
 \label{Alg:simple_presolve}
\end{algorithm}

The complexity of Alg~\ref{Alg:simple_presolve} is $O(\mbox{NNZ})$, where $\mbox{NNZ}$ is the number of nonzero elements in $A$. We implement this procedure solely in serial since it executes very quickly in practice.

\subsection{Parallel Maximal Clique Detection from Knapsack Constraint}
\label{sec:kcgen}
Following PBC generation, we proceed to detect maximal cliques from conflicting knapsack constraints with the Clique Detection method of Brito and Santos \cite[Algorithm 1]{brito2021preprocessing}. We modify the method, described as (serial) Alg~\ref{Alg:clique_detect}: namely, instead of returning all detected maximal cliques together, we separately return the first detected maximal clique (see line 7 in Alg~\ref{Alg:clique_detect}) as $\mathcal{S}_{org}$ (the \emph{original maximal clique}) and all other cliques in $\mathcal{S}_{other}$ (\emph{other maximal cliques}). This is used for more intensive clique extension and merging applied specifically to the original cliques, described in Section~\ref{sec:parext}. Furthermore, this method can produce a set $S_{other}$ that takes a quadratic amount of memory for storing all completed cliques; thus we set a size limit on knapsacks to avoid out of memory (see Sec.~\ref{sec:complim}). We note that the data structure proposed by Brito and Santos \cite{brito2021preprocessing} involves shared lists that seem difficult to parallelize, hence our alternative choice. Alg~\ref{Alg:clique_detect} for $\mathcal{S}_{ck}$ is, in turn, called in parallel via Alg~\ref{Alg:parallel_clique_detect}.

\begin{algorithm}[!htbp]
\SetAlgoLined
\LinesNumbered
\SetKwRepeat{Do}{do}{while}
\SetKwInput{Input}{Input}
\SetKwInput{Output}{Output}
\Input{Knapsack Constraint $\sum_{j\in\mathcal{B}}a_jx_j \leq b$}
\Output{$\mathcal{S}_{org}$, $\mathcal{S}_{other}$}
    Sort index set $\mathcal{B} = \{j_1, ...,j_n\}$ by non-decreasing coefficient value $a_{j_1} \leq ...\leq a_{j_n}$\;
    \If{$a_{j_{n-1}} + a_{j_n} \leq b$}{\textbf{return} $\emptyset, \emptyset$}
    Set $\mathcal{S}_{org}$, $\mathcal{S}_{other} := \emptyset$\;
    Find the smallest $\phi$ such that $a_{j_\phi} + a_{j_{\phi+1}} > b$\;
    Set $\mathcal{S}_{org} := \{x_{j_\phi},...,x_{j_n}\}$\;
    \For{$i = \phi-1 : 1$}{
        Find the smallest $\sigma$ such that $a_{j_i} + a_{j_\sigma} > b$\;
        \eIf{$\sigma$ exists}{
            Set $\mathcal{S}_{other} := \mathcal{S}_{other} \cup \{x_{j_i}, x_{j_\sigma},...,x_{j_n}\}$\;
        }{
            \textbf{Break}\;
        }
    }
 \textbf{return} $\mathcal{S}_{org}$, $\mathcal{S}_{other}$.
 \caption{Clique Detection}
 \label{Alg:clique_detect}
\end{algorithm}

Furthermore, in line 2 of Alg~\ref{Alg:parallel_clique_detect}, we randomly shuffle and then partition $\mathcal{S}_{ck}$ into subsets of equal cardinality (modulo the last processor) from Alg~\ref{Alg:simple_presolve} as a simple and fast heuristic for load balancing. Moreover, this can be efficiently parallelized \cite{bacher2015mergeshuffle,DBLP2015}. The same trick is used in Alg~\ref{Alg:parallel_cg} and Alg~\ref{Alg:parallel_extend_clique}. Note that the general problem of optimal load balancing---dividing up tasks with known computational load as evenly as possible across $k$ cores---is an NP-hard partitioning problem \cite{chopra1993partition,devine2006partitioning}.  The shuffling heuristic is justified both by average-case analysis as well as computational experiments indicating high parallel efficiency on hard instances. 

\begin{algorithm}[!htbp]
\SetAlgoLined
\LinesNumbered
\SetKwRepeat{Do}{do}{while}
\SetKwInput{Input}{Input}
\SetKwInput{Output}{Output}
\Input{$\mathcal{S}_{ck}$, $k$ threads}
\Output{$\mathcal{C}_{org}$, $\mathcal{C}_{other}$}
 Set $\mathcal{C}_{org}$, $\mathcal{C}_{other}:= \emptyset, \emptyset$\;
 Randomly shuffle $\mathcal{S}_{ck}$ and partition evenly by cardinality into subsets $\mathcal{S}_{ck}^{1},...,\mathcal{S}_{ck}^{k}$\;
 \For(parallel){$i \in \{1,...,k\}$}{
    Set $\mathcal{C}_{org}^i$, $\mathcal{C}_{other}^i:= \emptyset$\;
    \For{$c \in \mathcal{S}_{ck}^{k}$}{
        Set $\mathcal{S}_{org}$, $\mathcal{S}_{other}$ from Clique Detection (Alg 2) for constraint $c$\;
        Set $\mathcal{C}_{org}^i := \mathcal{C}_{org}^i \cup \{\mathcal{S}_{org}\}$, $\mathcal{C}_{other}^i := \mathcal{C}_{other}^i \cup \{\mathcal{S}_{other}\}$\;
    }
    Set $\mathcal{C}_{org} := \mathcal{C}_{org} \cup \mathcal{C}_{org}^i$, $\mathcal{C}_{other} := \mathcal{C}_{other} \cup \mathcal{C}_{other}^i$\;
 }
 \textbf{return} $\mathcal{C}_{org}$, $\mathcal{C}_{other}$. 
 \caption{Parallel Clique Detection}
 \label{Alg:parallel_clique_detect}
\end{algorithm}

Algorithm~\ref{Alg:clique_detect} has a worst-case runtime complexity of $O(n^2)$\cite[Page 4, Paragraph 2-3]{brito2021preprocessing}, where $n$ is the number of elements in the constraint. For Algorithm~\ref{Alg:parallel_clique_detect}, suppose we have $|\mathcal{S}_{ck}|:=m$ knapsacks and $k$ threads. The random shuffle of $m$ knapsacks can be completed within $O(m/k)$ \cite{bacher2015mergeshuffle,DBLP2015}. In lines 3-10, the workload of each thread is $O(mn^2/k)$, which is the worst-case runtime complexity of Algorithm~\ref{Alg:parallel_clique_detect}.

\subsection{Parallel Conflict Graph Construction}
\label{sec:cgcons}
Following clique detection from PBC-derived conflicting knapsack constraints, we proceed to CG construction. Suppose there are $n_{\mathcal{B}}$ binary variables in $\mathcal{B}$, and so  $n_{\mathcal{B}}$ complementary variables. We represent the CG with a sparse matrix $G\in \{0,1\}^{2n_{\mathcal{B}}\times 2n_{\mathcal{B}}}$, with the first set of rows $j\in \{1,...,n_{\mathcal{B}}\}$ representing the original variables $x_j$ and the next set $j\in \{n_{\mathcal{B}}+1,...,2n_{\mathcal{B}}\}$ representing the complements $\bar{x}_{j-n_{\mathcal{B}}}$. We build $G$ by Alg~\ref{Alg:parallel_cg} in parallel. 

We choose a sparse adjacency matrix here instead of a clique table in order to enable faster clique extension (see Section~\ref{sec:parext} for details) by avoiding redundant computation. A clique table, however, is substantially more efficient with memory, so as a workaround to the sparse data structure we set limits on the number of nonzeros considered (see Section~\ref{sec:complim} for details). Moreover, on cliques that we choose not to extend (namely $C_{other}$), we adopt the more memory-efficient data structure described in Section 2.3 of \cite{brito2021preprocessing}. 
 
The initialization of the clique set $\mathcal{C}$ in line 1 of Algorithm~\ref{Alg:cg_build} includes: the trivial pairwise conflicts between variables and their complements (i.e. $G_{j, j+n_{\mathcal{B}}} = 1$ for all $j=1,...,n_{\mathcal{B}}$); the cliques of $\mathcal{S}_{osp}$ from set packing constraints obtained by Alg~\ref{Alg:simple_presolve}; and cliques extracted in Alg~\ref{Alg:parallel_clique_detect}.

In line 2, we randomly shuffle and partition the clique set as a heuristic for parallel load balancing

In line 5, each $G(i)$ is updated in parallel using the clique subset $\mathcal{C}^i$ via Alg~\ref{Alg:cg_build}.

In lines 7-12, we reduce the individual $G(i)$ to the global $G$ via binary combination (see Fig~\ref{fig:parallel_cg}). In line 9, the \texttt{OR} logical merge is applied to every element in CGs, e.g., $G(i)_{j_1, j_2} := \max\{G(i-2^{d-1})_{j_1, j_2}, G(i)_{j_1,j_2}\}$.

\begin{algorithm}[!htbp]
\SetAlgoLined
\LinesNumbered
\SetKwRepeat{Do}{do}{while}
\SetKwInput{Input}{Input}
\SetKwInput{Output}{Output}
\Input{Clique set $\mathcal{C}$}
\Output{G}
    Initialize $G$ as an all-zero matrix with a size $2n_{\mathcal{B}}\times 2n_{\mathcal{B}}$\;
    \For{$Q \in \mathcal{C}$}{
        \For{$\forall (x_i,x_j)\ \mbox{s.t.}\ x_i, x_j\in Q$}{
            Set $G_{ij} := G_{ij}+1$\;
        }
    }
 \textbf{return} G.
 \caption{Conflict Graph Construction}
 \label{Alg:cg_build}
\end{algorithm}

\begin{algorithm}[!htbp]
\SetAlgoLined
\LinesNumbered
\SetKwRepeat{Do}{do}{while}
\SetKwInput{Input}{Input}
\SetKwInput{Output}{Output}
\Input{clique set $\mathcal{C}$, $k$ threads, dimension $n_{\mathcal{B}}$}
\Output{$G$}
 Initialize $G$ as an all-zero matrix with a size $2n_{\mathcal{B}}\times 2n_{\mathcal{B}}$\;\;
 Randomly shuffle $\mathcal{C}$ and partition evenly by cardinality into subsets $\mathcal{C}^{1},...,\mathcal{C}^{k}$\;
 \For(parallel){$i \in \{1,...,k\}$}{
    Call Algorithm~(\ref{Alg:cg_build}) with $C^{i}$ to get $G(i)$\;
 }

 \For{$d \in \{1, ..., \lceil \log_2(k) \rceil\}$}{
    \For(parallel){$i \in \{2^d, 2*2^{d}, 3*2^{d}, ..., \lfloor k/(2^d)\rfloor*2^d\}$}{
        Set $G(i) :=G(i-2^{d-1})\ \texttt{OR}\ G(i)$\;
    }
 }
 Set $G :=  G \ \texttt{OR}\ G(1)$\;
 \textbf{return} $G$. 
 \caption{Parallel Conflict Graph Construction}
 \label{Alg:parallel_cg}
\end{algorithm}

\begin{figure}[!htbp]
\centering
    \begin{tikzpicture}[scale=.75,auto=left]
    \path (5.5,3) node(x) {\color{red}Update CGs};
    \path (5.5,0) node(y) {\color{red}Combine CGs};
    \path (-6.2, 3) node(z) {\color{red} Thread Number};
        \node[circle,draw] (a) at (-4,3) {1}; 
        \node[circle,draw] (b) at (-3,3) {2}; 
        \node[circle,draw] (c) at (-2,3) {3}; 
        \node[circle,draw] (d) at (-1,3) {4}; 
        \node[circle,draw] (e) at (0,3){5}; 
        \node[circle,draw] (f) at (1,3) {6}; 
        \node[circle,draw] (g) at (2,3) {7};
        \node[circle,draw] (h) at (3,3) {8};
        \node[circle,draw] (i) at (-3,1.5) {2};
        \node[circle,draw] (j) at (-1,1.5) {4};
        \node[circle,draw] (k) at (1,1.5) {6};
        \node[circle,draw] (l) at (3,1.5) {8};
        \node[circle,draw] (m) at (-1,0) {4};
        \node[circle,draw] (n) at (3,0) {8};
        \node[circle,draw] (o) at (3,-1.5) {8};
        \draw[-,black] (a) -- (i);
        \draw[-,black] (b) -- (i);
        \draw[-,black] (c) -- (j);
        \draw[-,black] (d) -- (j);
        \draw[-,black] (e) -- (k);
        \draw[-,black] (f) -- (k);
        \draw[-,black] (g) -- (l);
        \draw[-,black] (h) -- (l);
        \draw[-,black] (i) -- (m);
        \draw[-,black] (j) -- (m);
        \draw[-,black] (k) -- (n);
        \draw[-,black] (l) -- (n);
        \draw[-,black] (m) -- (o);
        \draw[-,black] (n) -- (o);
        \draw[->, blue] (h) -- (x);
        \draw[->, blue] (l) -| (y);
        \draw[->, blue] (n) -- (y);
        \draw[->, blue] (o) -| (y);
        \draw[->, blue] (z) -- (a);
    \end{tikzpicture}
    \caption{Parallel Structure for Algorithm~\ref{Alg:parallel_cg}}
    \label{fig:parallel_cg}
\end{figure}

Suppose that the clique set $C$ is a randomly generated set of cliques $Q_i \in C$. Further, suppose that the cliques $Q_i$ are randomly generated such that each variable $x_j$ is a member of $Q_i$ with a probability $p$ that follows a Bernoulli distribution. Then the average runtime complexity of Alg~\ref{Alg:parallel_cg} is $O(mn_{\mathcal{B}}^2p^2\log(mn_{\mathcal{B}}p)/k + \log_2 k \cdot n_{\mathcal{B}}^2)$ (see Proposition~\ref{prop:parallel_cg_complexity_average} in Appendix~\ref{app:theorems}). Furthermore, the worst case runtime complexity of Alg~\ref{Alg:parallel_cg} is $O(n^2 + \log_2 k \cdot n_{\mathcal{B}}^2)$, where $n$ is number of variables of Problem~(\ref{eq:MIP}) (see Corollary~\ref{cor:parallel_cg_complexity_worstcase} in Appendix~\ref{app:theorems}).

\subsection{Parallel Clique Extension and Merging}
\label{sec:parext}
After generating the conflict graph, we apply clique strengthening  \cite{achterberg2007conflict,achterberg2020presolve,brito2021preprocessing}---in particular, we modify the Clique Extension of Brito and Santos \cite{brito2021preprocessing}: a greedy algorithm used to generate one strengthened clique based on the original clique and the CG (see \cite[Algorithm 2]{brito2021preprocessing}). Our modification (see Alg~\ref{Alg:extend_clique}) involves (potentially) multiple clique strengthenings, and the increase in computational load is made practical by parallelization.

In line 1 of Alg~\ref{Alg:extend_clique}, we generate a list $L$ including all variables $u$ that do not belong to the clique, but for which $u$ conflicts with all variables in the clique. In lines 19-20, we isolate a largest extended clique $\mbox{Clq}^\prime$ from the clique set $\mathcal{C}$. We distinguish this longest extended clique from all other cliques in order to perform cut management, described in Section~\ref{sec:cutmgmt}.


\begin{algorithm}[!htbp]
\SetAlgoLined
\LinesNumbered
\SetKwRepeat{Do}{do}{while}
\SetKwInput{Input}{Input}
\SetKwInput{Output}{Output}
\Input{clique $\mbox{Clq}$, CG}
\Output{longest clique $\mbox{Clq}^{\prime}$, extended clique set $\mathcal{C}$}
 Set $L := \{u\not\in \mbox{Clq}: \mathrm{CG}_{uv} = 1, \forall v\in \mbox{Clq}\}$\;
 \If{$L = \emptyset$}{\textbf{return} Clq, $\emptyset$.}
 Set $\mathcal{C}:= \emptyset$\;
 \For{$u\in L$}{
    \eIf{$\mathcal{C} = \emptyset$}{
        Set $\mathcal{C} := \{\{u\}\}$\;
    }{
        \For{$Q \in \mathcal{C}$}{
            \If{$\mathrm{CG}_{uw}>0, \forall w\in Q$}{
                Set $Q := Q\cup\{u\}$\;
            }    
        }
        \If{$u$ is not added to any $Q\in\mathcal{C}$}{
            Set $\mathcal{C}:=\mathcal{C}\cup\{\{u\}\}$\;
        }
    }
 }
 Find a longest $Q^\prime\in \mathcal{C}$ and remove it from $\mathcal{C}$\;
 Set $\mbox{Clq}^\prime := \mbox{Clq} \cup Q^\prime$\;
 \For{$Q \in \mathcal{C}$}{
     Set $Q := \mbox{Clq} \cup Q$\;    
 }
 \textbf{return} $\mbox{Clq}^\prime, \mathcal{C}$. 
 \caption{Clique Extension}
 \label{Alg:extend_clique}
\end{algorithm}

Strengthening each clique is an independent procedure, and so we propose to apply strengthening in parallel as Alg~\ref{Alg:parallel_extend_clique}. 

\begin{algorithm}[!htbp]
\SetAlgoLined
\LinesNumbered
\SetKwRepeat{Do}{do}{while}
\SetKwInput{Input}{Input}
\SetKwInput{Output}{Output}
\Input{clique set $\mathcal{C}$, CG, $k$ threads}
\Output{longest extended clique set $\mathcal{C}^{long}$, other extended clique set $\mathcal{C}^{other}$}
 Set $\mathcal{C}^{long}, \mathcal{C}^{other} := \emptyset$\;
 Randomly shuffle $\mathcal{C}$ and partition evenly by cardinality into subsets $\mathcal{C}_1,...,\mathcal{C}_k$\;
 \For(parallel){$i = \{1,...,k\}$}{
    Set $\mathcal{C}^{long}_{i}, \mathcal{C}^{other}_{i} := \emptyset$\;
    \For{$Q \in \mathcal{C}_i$}{
        Set $\mbox{Clq}^\prime, \mathcal{C}^\prime := \texttt{Clique Extension}(Q, CG)$\;
        Set $\mathcal{C}^{long}_{i} := \mathcal{C}^{long}_{i} \cup \{\mbox{Clq}^\prime\}$\;
        Set $\mathcal{C}^{other}_{i} := \mathcal{C}^{other}_{i} \cup \mathcal{C}^\prime$\;
    }
    Set $\mathcal{C}^{long} := \mathcal{C}^{long}\cup \mathcal{C}^{long}_i$\;
    Set $\mathcal{C}^{other} := \mathcal{C}^{other}\cup \mathcal{C}^{other}_i$\;
 }
 \textbf{return} $\mathcal{C}^{long}, \mathcal{C}^{other}$. 
 \caption{Parallel Clique Extension}
 \label{Alg:parallel_extend_clique}
\end{algorithm}

Complexity for Alg~\ref{Alg:extend_clique} is $O(n_{\mathcal{B}}^2)$ (see Lemma~\ref{le:clique_extend_complexity} in Appendix~\ref{app:theorems}) for one clique; thus $O(mn_{\mathcal{B}}^2)$ for $m$ cliques. Alg~\ref{Alg:parallel_extend_clique} has complexity $O(mn_{\mathcal{B}}^2/k)$ (see Proposition~\ref{prop:parallel_clique_extend_complexity} in Appendix~\ref{app:theorems}). 

After obtaining extended cliques from Alg~\ref{Alg:parallel_extend_clique}, we check for domination between extended cliques and remove dominated cliques. For a clique defined as a set pack row/constraint (Equation~(\ref{eq:setpack})), we will only store the index set $\mathcal{S}$; thus, to check the domination between two cliques $\mathcal{S}_1$ and $\mathcal{S}_2$, we only need to check whether 
$\mathcal{S}_1\subseteq \mathcal{S}_2$ or $\mathcal{S}_2\subseteq \mathcal{S}_1$, which can be performed in $O(n_{\mathcal{B}})$ due to sparse data structures.


Given $m$ cliques, the domination checking process is in $O(m^2n_{\mathcal{B}})$. The parallel version of this,  Alg~\ref{Alg:parallel_clique_merge}, has complexity $O(m^2n/k)$. 

\begin{algorithm}[!htbp]
\SetAlgoLined
\LinesNumbered
\SetKwRepeat{Do}{do}{while}
\SetKwInput{Input}{Input}
\SetKwInput{Output}{Output}
\Input{clique set $\mathcal{C}$, $k$ threads}
\Output{modified $\mathcal{C}^\prime$}
 \For(parallel){$\forall (Q_i,Q_j)\ \mathrm{ s.t. }\ Q_i,Q_j\in \mathcal{C}$}{
    \uIf{$Q_i\ \mathrm{ dominates }\ Q_j$}{Remove $Q_j$ from $\mathcal{C}$\;}\ElseIf{$Q_j\ \mathrm{ dominates }\ Q_i$}{Remove $Q_i$ from $\mathcal{C}$\;}
 }
 \textbf{return} modified $\mathcal{C}$. 
 \caption{Parallel Clique Merging}
 \label{Alg:parallel_clique_merge}
\end{algorithm}

\subsection{Cut Triage}
\label{sec:cutmgmt}
The previous subsections describe subroutines that are run in presolve.  Subsequently, our cut triage heuristic sends key inequalities from the longest cliques directly to the original formulation (at root node), and all other cuts are applied dynamically throughout the search tree via lazy user cuts.  The generated cuts are categorized as follows:

In Alg~\ref{Alg:simple_presolve}, three constraint sets are generated: $\mathcal{S}_{osp}, \mathcal{S}_{isp}, \mathcal{S}_{ck}$. Subsequently, cliques are extracted from $\mathcal{S}_{ck}$ using Alg~\ref{Alg:parallel_clique_detect}, yielding $\mathcal{C}_{org}$ and $\mathcal{C}_{other}$. Applying Alg~\ref{Alg:parallel_cg} on the four sets $\mathcal{S}_{osp}, \mathcal{S}_{isp}, \mathcal{C}_{org}, \mathcal{C}_{other}$, the CG is constructed. Cliques are then strengthened/extended from the three sets $\mathcal{S}_{osp}, \mathcal{S}_{isp}, \mathcal{C}_{org}$ with parallel Alg~\ref{Alg:parallel_extend_clique}. As a result, this process yields $8$ clique sets: $\mathcal{C}^{long}_{osp}, \mathcal{C}^{other}_{osp}$ (from $\mathcal{S}_{osp}$), $\mathcal{C}^{long}_{isp}, \mathcal{C}^{other}_{isp}$ (from $\mathcal{S}_{isp}$), $\mathcal{C}^{long}_{org}, \mathcal{C}^{other}_{org}$ (from $\mathcal{C}_{org}$), and $\mathcal{C}_{other}^{long}$, $\mathcal{C}^{long}_{other}$ (from $\mathcal{C}_{other}$). However, the number of cliques from these sets can be quite large. Thus, adding all inequalities directly as constraints to the formulation is impractical both due to resource limitations on linear programming solves as well as potential numerical issues from excessive cuts. 

Our cut triage procedure applies a simple heuristic: adding inequalities as standard constraints to the original formulation until the total number of constraints has doubled. All additional inequalities beyond doubling are added as lazy user cuts.  Inequalities are considered in the following sequence:

First, we replace $\mathcal{S}_{org}$ in the MIP formulation (deleting it in line 5 of Alg~\ref{Alg:simple_presolve}) with the strengthened constraints from $\mathcal{C}^{long}_{osp}$. 

Second, all cliques from $\mathcal{C}_{osp}^{other}$, $\mathcal{C}_{isp}^{other}$, $\mathcal{C}_{org}^{other}$, and $\mathcal{C}_{other}^{other}$ are incorporated via \emph{user cuts}, which are placed in the cut pool and may be added to the model at any node in the branch-and-cut search tree to cut off relaxation solutions (see e.g. \cite[Page 769, Lazy Attribute]{gurobi}).

Third, for cliques from  $\mathcal{C}_{org}^{long}$, $\mathcal{C}_{isp}^{long}$ and $\mathcal{C}_{other}^{long}$, we either add them to the original formulation as constraints or else give them to the solver as \emph{user cuts} according to the following criteria:
\begin{enumerate}
    \item Let $\|C\|$ indicate the number of nonzero terms in $\mathcal{C}$, and $\mbox{NNZ}$ represents the number of nonzero terms in the constraints matrix, i.e. $A$, of MIP. Let $Clq\_\mbox{NNZ} := 0$.
    \item If $\|\mathcal{C}_{org}^{long}\|/\mbox{NNZ} \leq 1$, we add cliques from $\mathcal{C}_{org}^{long}$ as constraints and $Clq\_\mbox{NNZ} := Clq\_\mbox{NNZ} + \|\mathcal{C}_{org}^{long}\|$.
    \item If $\|\mathcal{C}_{isp}^{long}\|/\mbox{NNZ} \leq 1$, we add cliques from $\mathcal{C}_{isp}^{long}$  as constraints  and $Clq\_\mbox{NNZ} := Clq\_\mbox{NNZ} + \|\mathcal{C}_{isp}^{long}\|$.
    \item If $(Clq\_{\mbox{NNZ}} + \|\mathcal{C}_{other}^{long}\|)/\mbox{NNZ} \leq 1$, we add cliques from $C_{other}^{long}$  as constraints.
\end{enumerate}

This heuristic attempts to prioritize some inequalities that may be more impactful by forcing inclusion into the formulation as constraints, leaving the remaining inequalities for Gurobi to manage in the search tree via user cuts.

\subsection{Limiting Parameters}
\label{sec:complim}
CG management can be time- and memory-consuming (see e.g. \cite[Page 491 Paragraph 2]{achterberg2020presolve}). For instance, in CG construction, the sparse adjacency matrix costs $O(n^2)$ memory for the clique with length $n$. However, we are not obliged to consider every possible variable, nor identify every clique, etc.; as such, to avoid excessive memory or time costs, we set the following limits.

For clique extension in Alg~\ref{Alg:clique_detect}, we ignore all knapsack constraints exceeding containing more than $5000$ variables. For CG construction in Alg~\ref{Alg:cg_build}, if the length of the clique $Q$ is greater than $1000$, we randomly select $1000$ elements from $Q$ to construct CG, and we at most process $2.5\times 10^7$ nonzero terms. For Clique Extension in Alg~\ref{Alg:parallel_extend_clique}, we let each thread process at most $1.25\times10^6$ nonzero terms and return generated extended cliques to the main thread. For Clique Merging in Alg~\ref{Alg:parallel_clique_merge}, if the number of cliques is more than $10^5$, we give up to conduct Alg~\ref{Alg:parallel_clique_merge}.

\section{Numerical Experiments}

All code and data can be found in our \href{https://github.com/foreverdyz/ParallelCliqueMerge}{repository}\footnote{https://github.com/foreverdyz/ParallelCliqueMerge}.

\subsection{Experimental Setup}

\subsubsection{Test Set}
Experiments are conducted on the benchmark set of the MIPLIB 2017 Collection \cite{miplib}, consisting of $240$ instances. Because our parallel presolve method focuses on set packing constraints and conflicting knapsack constraints, we remove $64$ instances where the Alg.~\ref{Alg:simple_presolve} procedure yields $2$ or fewer such constraints. The results below are run on the remaining $176$ cases. 

\subsubsection{Software and Hardware}
All algorithms are implemented in Julia 1.10.2 \cite{bezanson2017julia} and a desktop running 64-bit Windows 11 with an AMD Ryzen Threadripper PRO 5975WX 32-Core CPU and 64 GB. This CPU has 32 physical cores and 64 logical processors. We solve both original MIPs and MIPs after applying the conflict graph management with Gurobi 11.0.0 \cite{gurobi} and an alpha (prototype) version of JuMP 1.28 \cite{jump}.

\subsubsection{Configurations}
Conflict Graph Management are run on $1$ (serial), $2$, $4$, $8$, $16$, $32$, $64$ threads.  For benchmarking experiments, we run Gurobi on a single core while allowing it to use up to 32 threads (default setting).

\subsubsection{Time Measurement}
Measured presolving time excludes overhead from the reading time of the input file as well as the runtime of Alg~\ref{Alg:simple_presolve} because Alg~\ref{Alg:simple_presolve} is already implemented in Gurobi \cite{achterberg2020presolve} and the focus of experiments is on the effects of novelties. Moreover, for the cases where $\mathcal{S}_{isp}$ and $\mathcal{S}_{ck}$ from Alg~\ref{Alg:simple_presolve} are empty sets, we do not deploy CG management since Gurobi has already implemented CG management for cliques from $\mathcal{S}_{osp}$---note that such cases are not included in our performance comparisons. Times are always given in seconds and represent wall-clock measurements. Furthermore, for all aggregations, we calculate the geometric mean---note that we apply a 1 second shift to runtime means, as is typical in MIP literature (e.g. \cite{miplib}).

\subsection{Parallel CG Performance} \label{sec:parallel_cg_exp}
In the following experiments we analyze the parallel efficiency of our parallel CG management algorithms (excluding branch-and-cut solver time). We set a time limit for the CG procedure of $120$ seconds. We exclude $78$ cases that can be trivially handled in less than $0.1$ seconds with single or multiple threads. On the remaining $99/176$ cases, parallel speed-ups (serial runtime divided by the parallel runtime) are presented in Fig.~\ref{fig:parallel_comp}, where the instances are separated based on the slowest CG runtime (with serial and different numbers of threads). $173/176$ cases take less than $30$ seconds in both serial and parallel, while the remaining $3$ cases achieve the time limit ($120$ seconds) in serial. 

\begin{figure}[!htbp]
\centering
\begin{tikzpicture}[scale=.7]
\begin{axis}[
    axis lines = left,
    xlabel={Number of Threads [$m$]},
    ylabel={Speed-up},
    xmin=1, xmax=64,
    ymin=1, ymax=16,
    xtick={1,2,4,8,16,32,64},
    ytick={1,2,4,8,16},
    legend pos=north west,
]
\addplot[
    color=YellowOrange,
    mark=*,
    style= ultra thick,
    ]
    coordinates {
    (1,1)(2, 1.5)(4,2.1)(8,2.7)(16,3.1)(32, 3.0)(64, 2.9)
    };
    ]
    \addplot[
    color=black,
    mark=pentagon*,
    style= ultra thick,
    ]
    coordinates {
    (1,1)(2, 1.57)(4, 2.7)(8,4.2)(16,5.9)(32, 7.8)(64, 8.1)
    };
    ]
    \addplot[
    color=cyan,
    mark=triangle*,
    style= ultra thick,
    ]
    coordinates {
    (1,1)(2, 1.7)(4, 3.4)(8, 5.5)(16,7.3)(32, 9.8)(64, 11.4)
    };
    ]
    \legend{Slowest Runtime $\geq 0.1$, Slowest Runtime $\geq 1$, Slowest Runtime $\geq 10$}
\end{axis}
\end{tikzpicture}
\caption{Speed-up vs. different number of threads, and different lines represent different ranges of serial runtimes (s).}
\label{fig:parallel_comp}
\end{figure}

\begin{figure}[!htbp]
\centering
\begin{tikzpicture}[scale=.7]
\begin{axis}[
    axis lines = left,
    xmode = log,
    log ticks with fixed point,
    xlabel={Runtime (seconds)},
    ylabel={Number of Presolved Cases},
    xmin=0, xmax=32,
    ymin=110, ymax=176,
    xtick={0.5, 1,2,4,8,16,32,64},
    legend pos=north west,
]
\addplot[
    color=Red,
    mark=square*,
    line width=0.8pt
    ]
    coordinates {
    (0.5, 110)(1.0, 129)(2.0, 139)(4.0, 151)(8.0, 160)(16.0, 171)(32.0, 173)
    };
    ]
    \addplot[
    color=blue,
    mark=o,
    line width=0.8pt
    ]
    coordinates {
    (0.5, 129)(1.0, 140)(2.0, 146)(4.0, 159)(8.0, 165)(16.0, 171)(32.0, 173)
    };
    ]
    \addplot[
    color=Orange,
    mark=diamond*,
    line width=0.8pt
    ]
    coordinates {
    (0.5, 128)(1.0, 138)(2.0, 149)(4.0, 159)(8.0, 164)(16.0, 171)(32.0, 176)
    };
    ]
    \addplot[
    color=Green,
    mark=star,
    line width=0.8pt
    ]
    coordinates {
    (0.5, 121)(1.0, 134)(2.0, 147)(4.0, 160)(8.0, 166)(16.0, 171)(32.0, 176)
    };
    ]
    \addplot[
    color=Cyan,
    mark=pentagon*,
    line width=0.8pt
    ]
    coordinates {
    (0.5, 115)(1.0, 134)(2.0, 146)(4.0, 160)(8.0, 166)(16.0, 172)(32.0, 176)
    };
    ]
    \legend{Serial, 2 Threads, 4 or 8 Threads, 16 or 32 Threads, 64 Threads}
\end{axis}
\end{tikzpicture}
\caption{Number of CG Processed Cases vs. different runtime limit, and different lines represent different numbers of threads.}
\label{fig:parallel_logplot}
\end{figure}

\subsection{Gurobi Performance}\label{sec:nocg}
In this section we consider the impact of our CG procedure on total solver time, using Gurobi as our benchmark. As aforementioned, we focus on 176 cases from MIPLIB for there are at least 3 packing/knapsack constraints for our CG procedure to work with.  We exclude from comparison in this subsection an additional 77 cases for which our method cannot derive additional inequalities beyond standard CG techniques, namely those instances for which $\mathcal{S}_{isp} \cup \mathcal{S}_{ck} = \emptyset$. For these 77 excluded cases our procedure will simply reduce to a standard (albeit parallelized) CG procedure that is already incorporated into most solvers' presolve routines (including Gurobi). Thus $99/173$ cases are considered for comparison purposes to analyze the impact of our CG routine on Gurobi. Each instance is run over 5 different random seeds in Gurobi to help mitigate the effects of solver variance, and so the results average over $99\times5 = 495$ instances. 

We present results in two parts: one with the $495$ instances excluding CG management/overhead time, and one with $495$ instances including CG management time and varying the number of threads used. Since our CG implementation will duplicate certain similar efforts in Gurobi's CG itself, we expect a fully-integrated implementation to have less overhead than what is reported; hence, the two times presented provide optimistic (without overhead), and conservative (with overhead) bounds on potential impact, respectively.

\subsubsection{Gurobi excluding CG management time}
Table~\ref{tab:nocg} compares Gurobi with and without our CG procedure: \emph{Runtime Comp.} is the total solve time of Gurobi+CG ignoring CG time divided by that of Gurobi alone; \emph{Nodes Comp.} is the same ratio in terms of number of search tree nodes; \emph{Faster} and \emph{Slower} are the number of solves for which the CG procedure led to $5\%$ faster or slower (respectively) runtimes. There are $51/495$ instances that Gurobi cannot solve within $3600$ seconds; $3/51$ of these can be solved by Gurobi+CG. Excluded from this table are $28/495$ instances that neither the original Gurobi nor Gurobi + CG can solve within the time limit, and we present computational details and comparison for the $28$ instances in Appendix~\ref{app:c}, Table~\ref{tab:gap_comp}. The column ``Bracket'' with ``$\geq x$ sec'' buckets instances according to the longest runtime between default Gurobi and Gurobi + CG.

\begin{table}[ht]
\centering
\caption{\label{tab:nocg}Performance comparison between default Gurobi and Gurobi + CG management.}
\setlength{\tabcolsep}{1.7mm}{
\begin{tabular}{lccccc}
\toprule[1pt]
Bracket &Instances  & Runtime Comp. &Nodes Comp. &Faster &Slower \\
\hline
All &467	&89.9\%	&66.8\%	&238	&155\\
$\geq 1$ sec &447	&89.6\%	&65.5\%	&227 &151\\
$\geq 10$ sec &373	&87.1\%	&59.9\%	&197	&126\\
$\geq 100$ sec &265	&84.2\%	&50.6\%	&143	&83\\
$\geq 1000$ sec &119 &81.5\%	&38.6\%	&57	&43\\
\bottomrule[1pt]
\end{tabular}}
\end{table}

\subsubsection{Gurobi including CG management time}
Table~\ref{tab:gurobi_cg} compares Gurobi with and without our parallel CG procedure: \emph{Runtime Comp.} is the total solve time of Gurobi+CG with CG time on different numbers of threads divided by that of Gurobi alone; \emph{Nodes Comp.} is the same ratio in terms of number of search tree nodes; \emph{Faster} and \emph{Slower} are the number of solves for which the CG procedure led to $5\%$ faster or slower (respectively) runtimes. Again, excluded from this table are $28/495$ instances  unsolved by either default Gurobi or Gurobi + CG within $3600$ seconds.

Table~\ref{tab:64cg} breaks out the results on 64 threads in terms of problem difficulty, in the same fashion as Table~\ref{tab:nocg}. The benefit of our parallel CG management procedure increases as problem difficulty increases. On such problems our CG procedure overhead is relatively negligible, as roughly the same reduction is reported in Table~\ref{tab:nocg}. This suggests that perhaps further gains on average speedup can be realized if problem difficulty can be reasonably estimated \emph{a priori}---for instance, by only activating parallel CG on sufficiently large instances.

Considering these total time results, our CG procedure has modest benefits in serial and demonstrates significant advantages as parallelization increases. 

\begin{table}[ht]
\centering
\caption{\label{tab:gurobi_cg}Performance comparison between default Gurobi and Gurobi + CG management on different threads.}
\setlength{\tabcolsep}{2mm}{
\begin{tabular}{lccccc}
\toprule[1pt]
Threads & Runtime Comp. &Faster &Slower \\
\hline
1 thread &98.5\%	&221 &182\\
2 threads &95.6\%	&225 &180\\
4 threads &94.5\%	&227 &180\\
8 threads &93.7\%	&227 &180\\
16 threads &93.0\% &228 &180\\
32 threads &92.6\% &230 &178 \\
64 threads &92.3\% &230 &178  \\
\bottomrule[1pt]
\end{tabular}}
\end{table}

\begin{table}[ht]
\centering
\caption{\label{tab:64cg}Performance comparison between default Gurobi and Gurobi + CG management on 64 threads.}
\setlength{\tabcolsep}{1.7mm}{
\begin{tabular}{lccccc}
\toprule[1pt]
Bracket &Instances  & Runtime Comp. &Nodes Comp. &Faster &Slower \\
\hline
All &467	&92.3\%	&66.8\%	&230	&178\\
$\geq 1$ sec &448	&92.0\% &65.5\%	&224 &165\\
$\geq 10$ sec &375	&88.4\% &59.9\%	&198	&128\\
$\geq 100$ sec &265	&84.4\% &50.6\%	&142	&85\\
$\geq 1000$ sec &119	&82.0\% &38.6\%	&57	&45\\
\bottomrule[1pt]
\end{tabular}}
\end{table}

\section{Conclusion}
We develop an efficient parallel CG scheme that generates a larger number of valid inequalities than could otherwise be practically found in serial. Computational experiments demonstrate that this approach yields substantial overall speedups to Gurobi primarily due to search tree node reduction from the larger pool of generated cuts and valid inequalities. Rather than accelerating an existing subroutine in solvers, we leverage parallelism and intensify the effort by modifying the serial subroutine to conduct more computations given the same time budget. This parallelization approach enables substantial parallel acceleration from personal computing setups, and could yield benefits on other subroutines used in branch-and-cut solvers.

\section*{Acknowledgements} This work was funded by the Office of Naval Research under grant N00014-23-1-2632. We sincerely thank Imre Pólik (FICO) and Gerald Gamrath (COPT) for their valuable insights and suggestions on handling sparse data in our implementation. We are also grateful to Benjamin Müller (COPT) for identifying errors related to reading MPS files in Julia JuMP. 
\bibliographystyle{splncs04}
\bibliography{references}

\newpage
\appendix

\section{Complexity Results}\label{app:theorems}

\begin{lemma}\label{le:cg_complexity_average}
    For a clique set $\mathcal{C}$, suppose that for each clique $Q\in\mathcal{C}$, whether a variable $x_j$ is in $Q$ is given by a Bernoulli distribution with independent probability $p$. Then Alg~\ref{Alg:cg_build} has an average runtime of $O(mn_{\mathcal{B}}^2p^2\log(mn_\mathcal{B}p))$, where $n_\mathcal{B}$ is the total number of binary variables, and $m$ is the number of cliques in $\mathcal{C}$.
\end{lemma}
\begin{proof}
    In line 2, there are $m$ cliques.
    
    For $Q\in\mathcal{C}$, let $n_{Q}$ be the number of binary variables in $Q$. Then $n_{Q}$ is in a binomial distribution $B(n_\mathcal{B}, p)$; thus $\mathbb{E}[n_{Q}] = n_\mathcal{B}p$ and $\mathbb{E}[n_Q^2] = n_\mathcal{B}p(1-p) + n_\mathcal{B}^2p^2$ \cite{binomial}. In line 3, there are $n_Q(n_Q-1)/2$ pairs of variables $(x_i,x_j)$. We have
    \begin{equation*}
        \begin{aligned}
            \mathbb{E}[n_Q(n_Q-1)] =& \mathbb{E}[n_Q^2-n_Q]\\ 
            =& \mathbb{E}[n_Q^2] - \mathbb{E}[n_Q]\\
            =& n_\mathcal{B}p(1-p) + n_\mathcal{B}^2p^2 - n_\mathcal{B}p\\
            =& (n_\mathcal{B}^2-n_\mathcal{B})p^2 \in O(n_{\mathcal{B}}^2p^2).
        \end{aligned}
    \end{equation*}
    Therefore, $\mathbb{E}[mn_Q(n_Q-1)/2] = m \mathbb{E}[n_Q(n_Q-1)/2] \in O(mn_{\mathcal{B}}^2p^2)$. The average time complexity to build a sparse matrix (CSC format adopted by Julia) with $\hat{n}$ nonzero terms is $O(\hat{n}\log(\hat{n}))$ \cite{davis2006direct}; therefore, Alg~\ref{Alg:cg_build} has an average runtime in $O(mn_{\mathcal{B}}^2p^2\log(mn_{\mathcal{B}}p))$.
    \qed
\end{proof}

\begin{proposition}\label{prop:parallel_cg_complexity_average}
     For a clique set $\mathcal{C}$, suppose that for each clique $Q\in\mathcal{C}$, whether a variable $x_j$ is in $Q$ is given by a Bernoulli distribution with the probability $p$. Then Alg~\ref{Alg:parallel_cg} has an average runtime of $O(mn_{\mathcal{B}}^2p^2\log(mn_{\mathcal{B}}p)/k +\log k\cdot n_{\mathcal{B}}^2)$, where $n_{\mathcal{B}}$ is the total number of binary variables, $m$ is the number of cliques, and $k$ is the number of threads.
\end{proposition}
\begin{proof}
    Randomly shuffling and partitioning in line 2 can be done in parallel in $O(m/k)$ \cite{bacher2015mergeshuffle,DBLP2015}.
    
    In lines 3-6, there is a parallel for loop involving $k$ threads. In each iteration, we call Alg~\ref{Alg:cg_build} to at most $\lceil\frac{m}{k}\rceil$ cliques. Due to Lemma~\ref{le:cg_complexity_average}, the average complexity of lines 3-6 is $O(mn_{\mathcal{B}}^2p^2/k\log(mn_{\mathcal{B}}p))$.
    
    In lines 7-11, there is a binary combination with $\lceil\log k\rceil$ iterations. In the $i$-th iteration, we perform $\lfloor k/2^i\rfloor$ \texttt{OR} operations in parallel. The runtime of each \texttt{OR} is in $O(n_\mathcal{B}^2)$; thus the complexity of lines 7-11 (the binary combination) is $O(\log k \cdot n_{\mathcal{B}}^2)$.

    Line 12 is another \texttt{OR} operation in which the complexity is $O(n_\mathcal{B}^2)$.

    Consequently, Alg~\ref{Alg:parallel_cg} has an average runtime of $O(mn_{\mathcal{B}}^2p^2\log(mn_{\mathcal{B}}p)/k +\log k\cdot n_{\mathcal{B}}^2)$.
    \qed
\end{proof}

\begin{corollary}\label{cor:parallel_cg_complexity_worstcase}
    Alg~\ref{Alg:parallel_cg} has a worst-case runtime of $O(n^2 +\log k\cdot n_{\mathcal{B}}^2)$, where $n_{\mathcal{B}}$ is the number of binary variables, $m$ is the number of cliques, and $k$ is the number of threads.
\end{corollary}
\begin{proof}
    The proof is similar that of Proposition~\ref{prop:parallel_cg_complexity_average}. The only difference is in lines 3-6. In the worst case, the time complexity to build CG in one thread is $O(n^2)$. Therefore, Alg~\ref{Alg:parallel_cg} has a worst case runtime of $O(n^2 +\log k\cdot n_{\mathcal{B}}^2)$.
\end{proof}

\begin{lemma}\label{le:clique_extend_complexity}
    Alg~\ref{Alg:extend_clique} has complexity $O(n_\mathcal{B}^2)$, where $n_\mathcal{B}$ is the total number of binary variables.
\end{lemma}
\begin{proof}
    A clique Clq has at most $n_\mathcal{B}$ variables. In line 1, checking whether $\mathrm{CG}_{uv} = 1$ for a given $u\not\in \mathrm{Clq}$ and all $v\in \mathrm{Clq}$ takes at most $O(n_\mathcal{B})$ operations. So line 1 has complexity $O(n_\mathcal{B}^2)$.

    Lines 6-19 check all elements from $L$, whose size is bounded by $n_\mathcal{B}$. For the $t$-th element in $L$, $\mathcal{C}$ contains $t-1$ elements; thus lines 10 and 11 at most check $t-1$ pairs of $CG_{uw}$.  Therefore, the complexity of lines 6-19 is $O(\sum_{t=1}^{|L|} (t-1)) \in O(n_\mathcal{B}^2)$.

    In line 19, the size of $\mathcal{C}$ is $O(n_\mathcal{B})$; thus line 19 can be completed in $O(n_\mathcal{B})$.

    As a result, Alg~\ref{Alg:extend_clique} has complexity $O(n_\mathcal{B}^2)$.
\end{proof}

\begin{proposition}\label{prop:parallel_clique_extend_complexity}
    Alg~\ref{Alg:parallel_extend_clique} has complexity $O(mn_\mathcal{B}^2/k)$, where $n_{\mathcal{B}}$ is the number of binary variables, $m$ is the number of cliques, and $k$ is the number of threads.
\end{proposition}
\begin{proof}
    Randomly shuffling and partitioning in line 2 is $O(m/k)$ \cite{bacher2015mergeshuffle,DBLP2015}.

    In lines 3-10, there is a parallel for loop. In each iteration, we perform Alg~\ref{Alg:parallel_extend_clique} for at most $\lceil m/k\rceil$ cliques. Due to Lemma~\ref{le:clique_extend_complexity}, the runtime complexity of lines 3-10 is $O(mn_\mathcal{B}^2/k)$.

    As a result, Alg~\ref{Alg:parallel_extend_clique} has complexity $O(mn_\mathcal{B}^2/k)$. \qed
\end{proof}

\section{Detailed Computational Results}

In Table~\ref{tab:presolve_full}, \emph{CG T.} refers to the runtime of the CG procedure with 64 threads, and \emph{Org. T.} and \emph{Re. T.} are the runtimes of the Gurobi solver for the original and reduced models, respectively. All runtimes are rounded to two decimal places. \emph{Total}, \emph{Added}, and \emph{Users} separately correspond to the total number of detected cliques, the number of cliques added to the formulation as constraints, and the number of cliques designated as user cuts. \emph{Org. N.} and \emph{Re. N.} indicate the number of nodes explored by Gurobi for the original and reduced models separately. 

\setlength{\tabcolsep}{0.3mm} 
\begin{longtable}{lcccccccc}
    \caption{Detailed information for each instance}
    \label{tab:presolve_full} \\
    \toprule
    Instance & CG T. & Total& Added & User &Org. T. &Re. T. &Org. N. &Re. N.\\
    \midrule
    \endfirsthead
    
    \multicolumn{9}{c}{{\bfseries -- Continued from previous page --}} \\
    \toprule
    Instance & CG T.& Total& Added & User &Org. T. &Re. T. &Org. N. &Re. N.\\
    \midrule
    \endhead
    
    \midrule \multicolumn{9}{r}{{Continued on next page}} \\
    \endfoot
    \bottomrule
    \endlastfoot
30n20b8 &1.35 &110 &110 &0 &2.19 &2.12 &248 &248\\
academictimetablesmall &0.25 &22540 &19380 &3160 &301.97 &1354.57 &10911 &3343\\
atlanta-ip &0.19 &9294 &2640 &6654 &1304.76 &621.3 &16545 &11066\\
bab2 &1.79 &98005 &98005 &0 &1687.84 &2003.42 &64485 &67674\\
bab6 &1.23 &87992 &87992 &0 &1486.01 &691.15 &3356 &2153\\
bnatt400 &0.01 &8037 &8037 &0 &216.65 &230.12 &39568 &46186\\
bnatt500 &0.01 &10077 &10077 &0 &730.62 &546.61 &128219 &83733\\
brazil3 &0.12 &33214 &32414 &800 &310.93 &85.88 &10931 &6778\\
co-100 &7.15 &1135727 &2105 &1133622 &2345.76 &1194.12 &615808 &3968\\
cryptanalysiskb128n5obj14 &0.74 &104954 &104954 &0 &2505.57 &936.36 &519 &423\\
cryptanalysiskb128n5obj16 &0.62 &104962 &104962 &0 &2690.44 &1289.46 &628 &533\\
csched007 &0.0 &141 &141 &0 &113.46 &112.12 &146458 &146458\\
csched008 &0.0 &132 &132 &0 &328.29 &321.97 &845503 &845503\\
decomp2 &0.03 &11907 &11907 &0 &0.72 &0.38 &1 &1\\
drayage-100-23 &0.05 &281 &281 &0 &0.19 &0.17 &1 &1\\
drayage-25-23 &0.05 &277 &277 &0 &3.07 &2.67 &1 &1\\
dws008-01 &0.02 &510 &510 &0 &2200.58 &3601.2 &433711 &510088\\
enlight\_hard &0.0 &6 &6 &0 &0.02 &0.01 &0 &0\\
fastxgemm-n2r6s0t2 &0.0 &32 &32 &0 &36.64 &35.64 &21778 &21778\\
fhnw-binpack4-4 &0.0 &294 &294 &0 &25.06 &12.95 &172294 &84466\\
fhnw-binpack4-48 &0.0 &1095 &1095 &0 &12.79 &10.15 &10358 &5607\\
germanrr &0.01 &135 &135 &0 &547.49 &1144.46 &45168 &79827\\
gfd-schedulen180f7d50m30k18 &2.17 &179841 &177751 &2090 &3602.38 &3672.25 &14222 &14973\\
gmu-35-40 &0.0 &424 &366 &58 &144.46 &31.12 &155624 &27699\\
gmu-35-50 &0.0 &472 &376 &96 &128.07 &57.46 &153143 &54190\\
highschool1-aigio &8.39 &568786 &568786 &0 &1266.26 &477.95 &3967 &1\\
hypothyroid-k1 &0.34 &213094 &213094 &0 &9.91 &5.25 &1 &1\\
irish-electricity &0.02 &9765 &9746 &19 &1088.96 &1462.29 &11818 &17077\\
istanbul-no-cutoff &0.0 &22 &22 &0 &42.54 &41.05 &128 &106\\
k1mushroom &4.28 &836701 &836701 &0 &179.24 &142.75 &1 &1\\
lectsched-5-obj &0.04 &1541 &1541 &0 &124.44 &80.31 &9356 &6006\\
mcsched &0.0 &2026 &2026 &0 &25.43 &21.28 &11414 &10226\\
mzzv11 &0.06 &6703 &6640 &63 &11.73 &11.19 &1 &1\\
mzzv42z &0.06 &7512 &7436 &76 &7.32 &7.4 &1 &1\\
n2seq36q &0.38 &65231 &1929 &63302 &2.29 &13.4 &1 &1\\
neos-2657525-crna &0.0 &91 &91 &0 &2881.5 &3600.35 &22605925 &41462019\\
neos-2746589-doon &0.44 &11760 &11760 &0 &31.44 &237.31 &1861 &1387\\
neos-3216931-puriri &0.01 &2149 &1746 &403 &99.09 &89.03 &3343 &3768\\
neos-3656078-kumeu &0.01 &1904 &1904 &0 &1820.0 &1973.35 &135707 &78828\\
neos-3988577-wolgan &0.09 &27300 &27300 &0 &240.95 &214.96 &3075 &3017\\
neos-4387871-tavua &0.01 &5985 &5985 &0 &685.45 &879.27 &58796 &60241\\
neos-4413714-turia &3.65 &190954 &190954 &0 &30.92 &34.42 &1 &1\\
neos-4532248-waihi &2.09 &86278 &86278 &0 &364.63 &182.52 &2603 &934\\
neos-4647030-tutaki &0.01 &4200 &4200 &0 &169.5 &70.47 &14943 &223\\
neos-4722843-widden &0.49 &76669 &76669 &0 &49.66 &33.45 &2736 &1938\\
neos-4738912-atrato &0.0 &133 &133 &0 &37.94 &24.46 &7761 &8332\\
neos-5104907-jarama &0.01 &4286 &4286 &0 &1053.06 &2231.28 &118 &149\\
neos-5195221-niemur &0.04 &16908 &16908 &0 &2564.73 &3601.43 &38526 &44275\\
neos-631710 &9.55 &264198 &169576 &94622 &1026.19 &3603.27 &16073 &9784\\
neos-662469 &0.1 &478 &478 &0 &25.38 &22.01 &95 &127\\
neos-787933 &4.29 &234612 &1765 &232847 &2.1 &4.46 &1 &1\\
neos-827175 &0.06 &5518 &5518 &0 &0.95 &0.6 &1 &1\\
neos-957323 &1.23 &4202 &2612 &1590 &2.7 &4.74 &1 &1\\
neos-960392 &0.7 &5502 &5162 &340 &3.89 &5.61 &1 &1\\
neos8 &0.58 &65916 &65916 &0 &0.67 &3.96 &1 &1\\
net12 &0.01 &827 &827 &0 &172.39 &137.79 &4448 &2833\\
netdiversion &1.3 &60006 &60006 &0 &267.15 &236.72 &38 &57\\
ns1116954 &0.05 &39951 &7377 &32574 &225.18 &331.34 &1 &1399\\
nu25-pr12 &0.01 &479 &479 &0 &3.8 &3.73 &33 &33\\
nursesched-medium-hint03 &0.65 &5647 &5438 &209 &1336.71 &3633.57 &12588 &297\\
nursesched-sprint02 &0.25 &2022 &1882 &140 &6.72 &5.25 &1 &1\\
peg-solitaire-a3 &0.02 &4340 &4340 &0 &121.59 &11.91 &9500 &1\\
physiciansched3-3 &0.66 &475050 &211892 &263158 &1172.26 &1306.84 &8033 &7499\\
physiciansched6-2 &1.2 &163087 &153414 &9673 &3.13 &16.09 &1 &1\\
piperout-08 &0.1 &4713 &4709 &4 &3.45 &3.73 &1 &1\\
piperout-27 &0.1 &6076 &6072 &4 &2.65 &2.92 &1 &1\\
proteindesign121hz512p9 &16.22 &233859 &290 &233569 &3045.68 &3616.62 &46226 &42954\\
proteindesign122trx11p8 &9.16 &196710 &243 &196467 &2330.22 &3606.67 &12042 &85960\\
radiationm18-12-05 &0.05 &12257 &12257 &0 &164.62 &325.51 &48150 &84611\\
rail01 &1.38 &37981 &37803 &178 &703.32 &323.1 &37 &55\\
rail02 &5.67 &78318 &77759 &559 &3600.99 &3601.82 &132 &62\\
rd-rplusc-21 &0.0 &165 &165 &0 &365.37 &385.98 &76944 &23200\\
rocI-4-11 &0.01 &3872 &3872 &0 &96.77 &28.45 &80351 &16115\\
rocII-5-11 &0.16 &7967 &1421 &6546 &770.8 &1466.95 &234409 &231217\\
roi2alpha3n4 &0.02 &13209 &13081 &128 &97.15 &52.26 &20492 &896\\
roi5alpha10n8 &0.92 &211612 &211047 &565 &3601.39 &1084.41 &156291 &4615\\
roll3000 &0.0 &798 &755 &43 &12.37 &12.93 &10306 &10390\\
s100 &11.19 &981 &981 &0 &3601.48 &3602.04 &12838 &14798\\
s250r10 &19.46 &10052 &10052 &0 &3603.53 &3603.74 &124971 &469956\\
satellites2-40 &0.12 &11871 &11871 &0 &369.54 &66.19 &1 &1\\
satellites2-60-fs &0.12 &11871 &11871 &0 &220.45 &285.75 &306 &744\\
sct2 &0.0 &1423 &1423 &0 &24.33 &35.7 &10545 &13915\\
sing326 &0.2 &14240 &14240 &0 &2300.28 &1838.24 &24734 &17923\\
sing44 &0.22 &14818 &14818 &0 &1791.49 &1378.66 &16124 &18182\\
splice1k1 &0.9 &875753 &875753 &0 &2702.34 &3600.13 &15107 &1\\
supportcase18 &0.11 &10015 &240 &9775 &2246.37 &127.22 &349328 &5661\\
supportcase22 &0.02 &4100 &4100 &0 &3601.56 &3600.98 &224 &239\\
supportcase33 &0.15 &210 &178 &32 &22.53 &21.11 &6497 &5511\\
supportcase6 &3.07 &279 &279 &0 &180.7 &179.69 &7177 &7177\\
supportcase7 &0.0 &460 &457 &3 &7.6 &7.39 &1 &1\\
swath1 &0.0 &21 &21 &0 &6.99 &5.88 &1413 &1413\\
swath3 &0.01 &21 &21 &0 &35.14 &50.15 &12575 &13394\\
tbfp-network &1.02 &74 &74 &0 &51.89 &52.08 &95 &95\\
traininstance2 &0.24 &2754 &2754 &0 &312.15 &112.74 &619146 &326573\\
traininstance6 &0.22 &2139 &2139 &0 &20.01 &18.76 &76860 &70634\\
triptim1 &0.09 &2073 &2060 &13 &12.33 &16.04 &1 &1\\
uccase9 &0.01 &1521 &1521 &0 &1590.91 &1536.24 &5658 &5658\\
uct-subprob &0.0 &562 &330 &232 &623.11 &372.95 &441205 &320599\\
wachplan &0.01 &672 &672 &0 &88.82 &83.13 &34835 &34835\\
\end{longtable}

\section{Computational Results for Exceed-Time-Limit Cases}\label{app:c}
In Table~\ref{tab:gap_comp}, \emph{Seed} denotes the random seed used for the experiments. \emph{Org. U.} and \emph{Red. U.} represent the upper bounds (i.e., the best feasible solutions found) from the original and reduced models, respectively. Likewise, \emph{Org. L.} and \emph{Red. L.} denote the corresponding lower bounds. The gap values, \emph{Org. Gap} and \emph{Red. Gap}, are computed as
\[Gap := \frac{\|U. - L.\|}{\max\{\|U.\|, \|L.\|\}},\]
where \emph{U.} and \emph{L.} corresponds to bounds from either the original or reduced model.

\begin{table}[!htbp]
\centering
\caption{\label{tab:gap_comp}Bounds comparison between default Gurobi and Gurobi + CG management on 64 threads for instances exceeding the time limit.}
\setlength{\tabcolsep}{1mm}{
\begin{tabular}{lccccccc}
\toprule[1pt]
Instance &Seed  &Org. U. &Org. L. &Red. U. &Red. L.  &Org. Gap &Red. Gap\\
\hline
gfd-schedulen180f7d50m30k18 &0	&-	&1	&- &1 &100\% &100\%\\
neos-2657525-crna &12	&1.81	&0	&1.81 &0 &100\% &100\%\\
neos-5195221-niemur &1234	&0.00384 &0.00247 &0.00384 &0.00156 &35.6\% &59.3\%\\
proteindesign121hz512p9 &12 &1473	&1454.7	&1473 &1469 &1.24\% &0.27\%\\
proteindesign121hz512p9 &123 &1473	&1453.6	&1473 &1461.2 &1.33\% &0.80\%\\
proteindesign122trx11p8 &1 &1450	&1742	&1747 &1744 &0.45\% &0.17\%\\
rail02 &0 & -198.42 &-203.33 &-194.56 &-204.52 &2.41\% &4.8\%\\  
rail02 &1 & -200.5 &-202.6 &-198.3 &-203.9 &1.04\% &2.75\%\\  
rail02 &12 & -198.42 &-203.7 &-182.0 &-204.52 &2.59\% &11.01\%\\  
rail02 &123 & -148.8 &-204.14 &-189.7 &-203.7 &27.11\% &6.87\%\\  
rail02 &1234 & -191.2 &-204.4 &-194.3 &-204.4 &6.46\% &4.94\%\\
s100 &0 &-0.169 &-0.17 &-0.169 &-0.17 &0.59\% &0.59\%\\ 
s100 &1 &-0.169 &-0.17 &-0.169 &-0.17 &0.59\% &0.59\%\\
s100 &12 &-0.169 &-0.17 &-0.169 &-0.17  &0.59\% &0.59\%\\ 
s100 &123 &-0.169 &-0.17 &-0.169 &-0.17  &0.59\% &0.59\%\\ 
s100 &1234 &-0.169 &-0.17 &-0.169 &-0.17  &0.59\% &0.59\%\\ 
s250 &0 &-0.171 &-0.172 &-0.171 &-0.172  &0.58\% &0.58\%\\ 
s250 &1 &-0.171 &-0.172 &-0.171 &-0.172  &0.58\% &0.58\%\\ 
s250 &12 &-0.171 &-0.172 &-0.171 &-0.172  &0.58\% &0.58\%\\ 
s250 &123 &-0.171 &-0.172 &-0.171 &-0.172  &0.58\% &0.58\%\\ 
s250 &1234 &-0.171 &-0.172 &-0.171 &-0.172  &0.58\% &0.58\%\\
supportcase22 &0 &- & 1.32 &- &1.58 &100\% &100\%\\
supportcase22 &1 &- & 0.964 &- &1.973 &100\% &100\%\\
supportcase22 &12 &- & 0.654 &- &1.936 &100\% &100\%\\
supportcase22 &123 &- & 1.33 &- &1.78 &100\% &100\%\\
supportcase22 &1234 &- & 0.511 &- &1.947 &100\% &100\%\\
splice1k1 &123 &-394 &-1490 &-394 &-1601 &73.56\% &75.39\%\\
splice1k1 &1234 &-394 &-1467 &-394 &-1609 &73.14\% &75.51\%\\
\bottomrule[1pt]
\end{tabular}}
\end{table}

From Table~\ref{tab:gap_comp}, the ratio of the geometric means of \emph{Org. Gap} to \emph{Red. Gap} is 96.58\%. This indicates that our cut generation (CG) management helps Gurobi achieve a slightly smaller primal-dual gap on the instances presented in the table.

\end{document}